\documentclass[12pt]{article}
\usepackage[top=1in, bottom=1in, left=0.5in, right=0.5in]{geometry}
\usepackage{hyperref}
\usepackage{amssymb}
\usepackage{amsmath}
\usepackage{amsthm}
\usepackage{mathrsfs}
\usepackage{graphicx}
\hypersetup{
    colorlinks=true,
    linkcolor=blue,
    filecolor=magenta,      
    urlcolor=cyan,
    pdftitle={Overleaf Example},
    pdfpagemode=FullScreen,
    }

\numberwithin{equation}{section}
\title{\textbf{\Large Algebraic analogues of results of Alladi-Johnson using the Chebotarev Density Theorem}}

\author{Sroyon Sengupta}
\date{}

\begin{document}

\maketitle
\textbf{Abstract:} \textit{{\small We aim to get an algebraic generalization of Alladi-Johnson's (A-J) work on Duality between Prime Factors and the Prime Number Theorem for Arithmetic Progressions - II, using the Chebotarev Density Theorem (CDT). It has been proved by A-J, that for all positive integers $k,\ell$ such that $1\leq \ell\leq k$ and $(\ell,k)=1$,}}
\begin{equation}
    \sum_{n\geq 2;\;p_1(n) \equiv \ell\;(mod\;k)}\frac{\mu(n)\omega(n)}{n} = 0, \nonumber
\end{equation}
\textit{{\small where $\mu(n)$ is the M\"obius function, $\omega(n)$ is the number of distinct prime factors of $n$, and $p_1(n)$ is the smallest 
prime factor of $n$. In our work here, we will prove the following result: If $C$ is a conjugacy class of the Galois group of some finite extension $K$ of $\mathbb{Q}$, then}}
\begin{equation}
    \sum_{ n \geq 2;\;\left[\frac{K/\mathbb{Q}}{p_1(n)}\right]=C} \frac{\mu(n)\omega(n)}{n} = 0. \nonumber
\end{equation}
\textit{{\small where $\left[\frac{K/\mathbb{Q}}{p_1(n)}\right]$ is the Artin symbol. When $K$ is a cyclotomic extension of $\mathbb{Q}$, this reduces to the exact case of A-J's result.}} 

\section{Introduction and Notation}
In 1977, Alladi [1] proved the following Duality Lemma: if f is a function defined on the primes, then:
\begin{eqnarray}
    \sum_{1<d|n} \mu(d)f(P_1(d)) = -f(p_1(n)), \nonumber\\
    \sum_{1<d|n}\mu(d)f(p_1(d)) = -f(P_1(n)), \nonumber
\end{eqnarray}
where $P_1(n)$ and $p_1(n)$ are the largest and smallest prime factors of $n$ respectively. He uses the second identity and properties of the M\"obius function in the later part of the paper to show that if $f$ is a bounded function on the primes such that 
\begin{eqnarray}
    \lim_{x \to \infty}\frac{1}{x} \sum_{2 \leq n \leq x}f(P_1(n)) = c, \nonumber
\end{eqnarray}
then
\begin{eqnarray}
    \sum_{n=2}^{\infty}\frac{\mu(n)f(p_1(n))}{n} = -c, 
\end{eqnarray}
and vice-versa. This is a surprising generalization of an important result by Landau, which states that
\begin{eqnarray}
    \sum_{n=1}^{\infty} \frac{\mu(n)}{n} = 0 
\end{eqnarray}
is (elementarily) equivalent to the Prime Number Theorem (PNT). To understand the generalization, we rewrite (1.2) as
\begin{eqnarray}
    \sum_{n=2}^{\infty} \frac{\mu(n)}{n}=-1. 
\end{eqnarray}
Note that from the general result in (1.1),  by taking $f(p)=1$ for all primes $p$, we get (1.3). \\ \\
It was also shown in [1] that the Prime Number Theorem in Arithmetic Progressions (PNTAP) implies that the sequence $P_1(n)$ of largest prime factors is uniformly distributed in the reduced residue classes modulo a positive integer $k$, which further means that if $f$ is chosen to be the characteristic function of primes in an arithmetic progression $\ell\;(mod\;k)$, then for each such $f$, the average of $f(P_1(n))$ exists and is $\frac{1}{\varphi(k)}$. Therefore, by (1.1) it follows that,  
\begin{eqnarray}
    \sum_{\substack{2 \leq n \\ p_1(n) \equiv \ell\;(mod\;k)}}^{\infty} \frac{\mu(n)}{n} = -\frac{1}{\varphi(k)} 
\end{eqnarray}
for all positive integers $k,\ell$ satisfying $(k,\ell)=1$. This is very surprising because it gives a way of slicing the convergent series in (1.3) into $\varphi(k)$ subseries, all converging to the same value! In [1], a strong quantitative version of (1.4) is also proved. \\ \\
Inspired by (1.4), Dawsey, in 2017 [4], proved similar density type results using the above duality, but in an algebraic setting of Galois extensions of the field of rational numbers $\mathbb{Q}$. Dawsey proved that if $C$ is a conjugacy class of the Galois group of some finite Galois extension $K$ of $\mathbb{Q}$, then
\begin{eqnarray}
    \lim_{x \to \infty}\sum_{\substack{2 \leq n \leq x \\ \left[\frac{K/\mathbb{Q}}{p_1(n)}\right] = C}} \frac{\mu(n)}{n} = -\frac{|C|}{|G|}.
\end{eqnarray}
In her paper, Dawsey uses the strong form of CDT given by Lagarias-Odlyzko [8], with effective bounds which are very similar to the ones in PNTAP, and hence is able to derive (1.5) in quantitative form. Dawsey considers the function $f$ to be characteristic function on the primes $p$ satisfying the condition that $\left[\frac{K/\mathbb{Q}}{p}\right] = C$. \textit{Theorem 2} in [Dawsey 4] states that the density of the integers satisfying the condition $\left[\frac{K/\mathbb{Q}}{P_1(n)}\right] = C$ is $\frac{|C|}{|G|}$, and therefore, by Alladi's general result, (1.5) will follow. However, since Dawsey uses the strong form of CDT by Lagarias-Odlyzko, she is able to get a quantitative form of (1.5), which (1.1) does not yield. \\ \\
In [1], Alladi also proved a general duality (see \textit{Theorem 1.3} below), which connects the smallest and the $k^{th}$ largest prime factors, and the dual which connects the largest to the $k^{th}$ smallest prime factor. Following is the version of one of his second order duality identities:
\begin{eqnarray}
    \sum_{d|n} \mu(d)(\omega(d)-1)f(p_1(d))= f(P_2(n)) 
\end{eqnarray}
Using (1.6), Alladi and his doctoral student Johnson [2] have recently proved the following: For integers $k,\ell$ such that $1 \leq \ell < k$ satisfying $(\ell,k)=1$,
\begin{eqnarray}
     \sum_{n\geq 2;\;p_1(n) \equiv \ell\;(mod\;k)}\frac{\mu(n)\omega(n)}{n} = 0 \nonumber
\end{eqnarray}
This is also a way of slicing the convergent series 
\begin{eqnarray}
    \sum_{n=2}^{\infty} \frac{\mu(n)\omega(n)}{n} = 0 \nonumber
\end{eqnarray}
into $\varphi(k)$ equal valued subseries. \\ \\
Motivated by Dawsey's work on the first order duality, in this paper, we aim to prove algebraic analogues of Alladi and Johnson's results in the setting of Galois extensions of the field of rational numbers, using the strong form of CDT. In the later part of [2], a more general result akin to (1.1) is stated. Our result could be viewed as a special case of that general result (see \S 5), but ours has been established in a quantitative form due to the use of Lagarias-Odlyzko CDT. \\ \\
To start with, let $K$ be a Galois extension of $\mathbb{Q}$ and $\mathcal{O}_K$ be the ring of integers of the field $K$. For an unramified prime $p$ and prime ideal $\mathfrak{p} \subset \mathcal{O}_K$ lying above $p$, let $\left[\frac{K/\mathbb{Q}}{p}\right]$ denote the Artin symbol. We will first recall the definition of Artin Symbol in generality for finite Galois extensions as done in Dawsey's paper [4]. Let $L/K$ is a finite Galois extension of number fields with corresponding number rings $\mathcal{O}_L$ and $\mathcal{O}_K$. For any non-zero prime ideal $\mathfrak{p}$ of $\mathcal{O}_K$, the ideal generated by $\mathfrak{p}$ in $\mathcal{O}_L$ can be uniquely decomposed into a product of 
distinct prime ideals $\mathfrak{p}_i$'s lying over $\mathfrak{p}$, i.e. 
\begin{equation}
    \mathfrak{p}\mathcal{O}_L = \prod_{i=1}^k \mathfrak{p}_i^{e_i} \nonumber
\end{equation}
Now, if $e_i=1$ for each $i$, we say that $\mathfrak{p}$ is unramified in $L$. This occurs for all but finitely many prime ideals $\mathfrak{p}$. We recall the definition of the absolute norm of a non-zero ideal $\mathcal{I}$ of a number ring $\mathcal{O}$ of some number field:
\begin{equation}
    Nm(\mathcal{I}) := [\mathcal{
    O}:\mathcal{I}]=\left|\frac{\mathcal{O}}{\mathcal{I}}\right| \nonumber
\end{equation}
We now have all the required machinery to define Artin Symbol, which is of  importance in our results. For every prime ideal $\mathfrak{p}_i$ in the decomposition of $\mathfrak{p}$, the Artin symbol $\left[\frac{L/K}{\mathfrak{p}_i}\right]$ for a particular factor $\mathfrak{p}_i$ is defined as the unique isomorphism $\sigma \in Gal(L/K)$ that satisfies the condition
\begin{equation}
    \sigma(\alpha) = \alpha^{Nm(\mathfrak{p})}\;(mod\;\mathfrak{p}_i) \nonumber
\end{equation}
for each $\alpha \in L$. Now, for every $\tau \in Gal(L/K)$, the ideals $\mathfrak{p}_i$'s are of course isomorphic to each other under each $\tau$ in some specific permutation depending on the choice of $\tau$. Further, we have that for such $\tau$,
\begin{equation}
    \left[\frac{L/K}{\tau(\mathfrak{p}_i)}\right] = \tau\left[\frac{L/K}{\mathfrak{p}_i}\right]\tau^{-1} \nonumber
\end{equation}
This implies that we can get a conjugacy class $C \subset Gal(L/K)$ which is associated to the prime ideal $\mathfrak{p}$ in $\mathcal{O}_L$. We define $\left[\frac{L/K}{\mathfrak{p}}\right]$, the Artin Symbol corresponding to $\mathfrak{p}$ as the associated conjugacy class $C$. Let us denote $\mathfrak{P}(K)$ to be the set of all prime ideals of $\mathcal{O}_K$. For a given conjugacy class $C \subset Gal(L/K)$, we define the set $\mathfrak{P}_C$ as follows:
\begin{eqnarray}
    \mathfrak{P}_C:=\left\{ \mathfrak{p} \in \mathfrak{P}(K):\;\mathfrak{p} \;\text{is unramified in $L$ and}\;  \left[\frac{L/K}{\mathfrak{p}_i}\right] = C  \right\} \nonumber
\end{eqnarray}
Let $\pi_{C}(x,L/K)$ denote the number of non-zero prime ideals $\mathfrak{p}$ in $\mathcal{O}_K$ that are unramified in $L$ with absolute norm bounded by $x$ whose Artin symbol is $C$. \\ \\
\textbf{Theorem 1.1:}(\textit{Chebotarev Density} [14]) \textit{Let $L/K$ be a finite Galois extension and $C$ be a conjugacy class of the Galois group $Gal(L/K)$. Then the natural density of $\mathfrak{P}_C$, defined by}
\begin{eqnarray}
    \lim_{x \rightarrow \infty}  \frac{\{\mathfrak{p} \in \mathfrak{P}_C:\;Nm(\mathfrak{p})\leq x\}}{\{\mathfrak{p}\in \mathfrak{P}(K):\; Nm(\mathfrak{p}) \leq x \}} \nonumber
\end{eqnarray}
\textit{exists and is equal to the ratio $\frac{|C|}{|Gal(L/K)|}$. More precisely, as $x \rightarrow \infty$}
\begin{eqnarray}
    \pi_C(x,L/K) = \frac{|C|}{|Gal(L/K)|}.\frac{x}{\log x} + o\left(\frac{x}{\log x}\right) \nonumber
\end{eqnarray}
A more precise formulation of the error has been given by Lagarias and Odlyzko [8] as follows: \\ \\
\textbf{Theorem 1.2:}(\textit{Lagarias-Odlyzko}) \textit{For sufficiently large $x \geq C$, where $C$ depends on both the absolute discriminant and the degree of extension of $L$ (over $\mathbb{Q}$), we have that}
\begin{eqnarray}
    \left|\pi_C(x,L/K)-\frac{|C|}{|Gal(L/K)|}Li(x)\right| <<x\exp\left\{-c_1\sqrt{\frac{\log x}{n_L}}\right\} \nonumber
\end{eqnarray}
\textit{for an appropriate constant $c_1$, where $n_L=[L:\mathbb{Q}]$ and $Li(x):=\int_2^x\frac{dt}{\log t}$.} \\
This bound for the error of the Chebotarev Density Theorem will be of significant importance to get our results. \\ \\
To demonstrate the problem, we now introduce a few important notations that we will use to prove our main results in this paper. We start with the general Alladi duality identities (Section 4 in [1]) that have paved the path for numerous results and is of major significance in the results that follow. Let $P_k(n)$ and $p_k(n)$ denote the $k^{th}$ largest prime factor of $n$ and the $k^{th}$ smallest prime factor of $n$ respectively, and if $\omega(n)<k$, then we set $P_k(n)=p_k(n)=1$. \\ \\
\textbf{Theorem 1.3:}(\textit{Alladi}) [1] \textit{If $f$ is an arithmetic function with $f(1)=0$, then the following four identities hold:}
\begin{eqnarray}
    \sum_{d|n} \mu(d)f(P_k(d))=(-1)^k {\omega(n)-1 \choose k-1}f(p_1(n)) \nonumber\\
    \sum_{d|n} \mu(d)f(p_k(d))=(-1)^k {\omega(n)-1 \choose k-1}f(P_1(n)) \nonumber
\end{eqnarray}
\textit{and}
\begin{eqnarray}
    \sum_{d|n} \mu(d){\omega(d)-1 \choose k-1}f(P_1(d))=(-1)^k f(p_k(n)) \nonumber\\
    \sum_{d|n} \mu(d){\omega(d)-1 \choose k-1}f(p_1(d))=(-1)^k f(P_k(n)) \nonumber
\end{eqnarray}
For our results, we will use the fourth duality identity with $k=2$, i.e. (1.6),
where, of course, $P_2(n)$ denotes the second largest prime factor of $n$. \\ \\
We utilize the following well-known function:
\begin{eqnarray}
    \Psi(x,y) = \sum_{n\leq x;\;P_1(n)\leq y} 1 \nonumber
\end{eqnarray}
i.e. $\Psi(x,y)$ counts the number of integers less or equal to $x$ whose largest prime factor is less or equal to $y$. We denote the quantity $\alpha$ to be the ratio $\frac{\log x}{\log y}$. We then have the following result by de Bruijn [6]:
\begin{eqnarray}
    \Psi(x,y) <<xe^{-c\alpha} \nonumber
\end{eqnarray}
uniformly for $2\leq y \leq x$, for some constant $c>0$. Tenenbaum [13] has shown that the above holds for $c=\frac{1}{2}$. de Bruijn [5] had further established the following significant uniform asymptotic estimate:
\begin{eqnarray}
    \Psi(x,y) \sim x\rho(\alpha) \nonumber
\end{eqnarray}
for $e^{(\log x)^{\frac{2}{5}}} \leq y \leq x$, where $\rho(\alpha)$ satisfies the following integro-difference equation
\begin{eqnarray}
    \rho(\alpha) = 1-\int_1^{\alpha} \frac{\rho(u-1)}{u}du\;\;\;\text{and} \;\;\;\rho(\alpha) = \exp\{-\alpha\log \alpha - \alpha \log \log \alpha +O(\alpha)\} \nonumber 
\end{eqnarray}
We shortly remark here that, the order of magnitude of $\Psi(x,y)$ is quite small in comparison to that of $x$ when $\alpha$ grows large.  \\ \\ 
We denote by $P_2(n)$, the second largest prime factor $n$. $P_2(n)$ could be defined in two different ways namely:
\begin{itemize}
    \item $P_2(n) = P_1(n/P_1(n))$; $P_2(n)=1$ when $\Omega(n)<2$,
    \item $P_2(n)$ is as the largest prime factor of $n$ strictly less than $P_1(n)$; $P_2(n)=1$ when $\omega(n) < 2$,
\end{itemize}
and so, there is an ambiguity. However, the following theorem in [2] shows that asymptotically, there is little difference between the above two definitions of $P_2(n)$. \\ \\
\textbf{Theorem 1.4:}
\textit{Let $N(x)$ denote the number of positive integers $n \leq x$ for which $P_1(n)$ repeats. Then}
\begin{eqnarray}
    N(x) << \frac{x}{\exp\{\frac{1}{2}\sqrt{\log x \log \log x}\}} \nonumber
\end{eqnarray}
This shows that numbers less than equal to $x$ with repeating largest prime factor is of the order $o(x)$, i.e. the set of integers for which the above definitions of $P_2(n)$ are ambiguous, is quite small as compared to $x$. In view of \textit{Theorem 1.4}, we understand that both the definitions of $P_2(n)$ are almost always same. Also, for all the theorems and their respective proofs that follow, we will only consider positive integers whose $P_1(n)$ does not repeat in the decomposition of the primes (i.e. exactly the set where the above two definitions of $P_2(n)$ are consistent).  \\ \\
We now need a quantitative version of the fact that for ``almost all" integers, $P_2(n)$ is large. This has been provided by Tenenbaum [12], on the size of $P_k(n)$, when $k \geq 2$. Tenenbaum's result is a stronger quantitative result on the joint distribution of the $P_k(n)$ for $k \geq 2$, of which, we only need the following moderation used in [2]: \\ \\
\textbf{Theorem 1.5:} \textit{For $2 \leq y \leq x$, we have}
\begin{eqnarray}
    \sum_{n \leq x;\;P_2(n) \leq y} 1 << \frac{x\log y}{\log x}. \nonumber
\end{eqnarray}
Alladi-Johnson established an elementary proof of \textit{Theorem 1.5}, when $y = e^{\sqrt{\log x}}$. We now have everything set up to move ahead with our main results. 
\section{Algebraic Analogues of Alladi-Johnson}
In this section, we aim to prove our main result. From here onwards, we specify that $K$ denotes a finite Galois extension of $\mathbb{Q}$, with Galois group $G = Gal(K/\mathbb{Q})$. Also, $C \subset G$ denotes an arbitrary, yet fixed, conjugacy class of the Galois group. We now state our first theorem, which we prove following the technique in [1] and [2]. \\ \\
\textbf{Theorem 2.1:} \textit{With the above notation in consideration, we have that}
\begin{eqnarray}
    N_2(x;K/\mathbb{Q},C)=\sum_{n \leq x;\;\left[\frac{K/\mathbb{Q}}{P_2(n)}\right]=C } 1 = \frac{|C|}{|G|}.x + O\left(\frac{x(\log \log x)^2}{\log x}\right) \nonumber
\end{eqnarray}
\textbf{Proof:} Firstly, we note that $P_2(n)=1$ if and only if $\omega(n)=1$ (since we are considering the second definition of $P_2(n)$). Now, positive integers up to $x$ satisfying $\omega(n)=1$ is given by
\begin{equation}
      \sum_{n \leq x;\;\omega(n)=1} 1 = \pi(x)+O(\sqrt{x}) = \frac{x}{\log x} + O\left(\frac{x}{\log^2x}\right) 
\end{equation}
We define the set $S_2(x,p)$ to be the set of integers $n \leq x$, such that $P_1(n)$ does not repeat and $P_2(n)=p>1$. We also define the set $S_2(x,C)$ to be the set of integers $n \leq x$, with $\left[\frac{K/\mathbb{Q}}{P_2(n)}\right]=C$ and $P_2(n)>1$. Therefore, denoting $P_2(n)=p$ and using Theorem 1.4 and (2.1), we have that
\begin{eqnarray}
     |S_2(x,C)| = \sum_{p \leq \sqrt{x};\;\left[\frac{K/\mathbb{Q}}{p}\right]=C} |S_2(x,p)| + O\left(\frac{x}{\log x}\right) 
\end{eqnarray}
From Theorem 1.5, we know that for any $y \leq x$
\begin{eqnarray}
    \sum_{p \leq y;\;\left[\frac{K/\mathbb{Q}}{p}\right]=C} |S_2(x,p)| \leq \sum_{p \leq y} |S_2(x,p)| << \frac{x\log y}{\log x}. \nonumber
\end{eqnarray}
Hence, using this in (2.2), we get that
\begin{eqnarray}
    |S_2(x,C)| = \sum_{y<p \leq \sqrt{x};\;\left[\frac{K/\mathbb{Q}}{p}\right]=C} |S_2(x,p)| + O\left(\frac{x\log y}{\log x}\right).
\end{eqnarray}
Let $N \in S_2(x,p)$ for some prime $p$. Then we can represent it $N=mpq$, where $P_1(N)=q$, $P_2(N)=p$ and $P_1(m)\leq p$. Here, $q$ is a prime such that $q>p$. Then of course, $m \leq \frac{x}{p^2}$, and hence, we can represent $|S_2(x,p)|$ as the following sum:
\begin{eqnarray}
    |S_2(x,p)| = \sum_{m \leq \frac{x}{p^2};\;P_1(m)\leq p} \;\;\sum_{p<q;\;mpq\leq x} 1 = \sum_{p<q\leq \frac{x}{p}}\;\;\sum_{m \leq \frac{x}{pq};\;P_1(m)\leq p} 1 = \sum_{p<q\leq\frac{x}{p}}\Psi\left(\frac{x}{pq},p\right) \nonumber
\end{eqnarray}
Therefore, from above we have that
\begin{eqnarray}
    \sum_{y<p \leq \sqrt{x};\;\left[\frac{K/\mathbb{Q}}{p}\right]=C} |S_2(x,p)| = \sum_{y<p\leq \sqrt{x};\;\left[\frac{K/\mathbb{Q}}{p}\right]=C} \;\; \sum_{p<q\leq \frac{x}{p}} \Psi\left(\frac{x}{pq},p\right) \hspace{6.3cm} \nonumber\\
    = \sum_{y < q \leq \sqrt{x}}\;\;\sum_{y<p<q;\;\left[\frac{K/\mathbb{Q}}{p}\right]=C} \Psi\left(\frac{x}{pq},p\right) + \sum_{\sqrt{x}<q\leq \frac{x}{y}}\;\;\sum_{y<p\leq \frac{x}{q};\;\left[\frac{K/\mathbb{Q}}{p}\right]=C} \Psi\left(\frac{x}{pq},p\right)   
\end{eqnarray}
We now replace the inner sums in the two separate terms on the right hand side of $(2.4)$ using the following respective integrals:
\begin{itemize}
    \item Replace $\sum_{y<p<q;\;\left[\frac{K/\mathbb{Q}}{p}\right]=C} \Psi\left(\frac{x}{pq},p\right)$ by $\frac{|C|}{|G|}\int_y^q\Psi\left(\frac{x}{tq},t\right)\frac{dt}{\log t}$
    \item Replace $\sum_{y<p<\frac{x}{q};\;\left[\frac{K/\mathbb{Q}}{p}\right]=C} \Psi\left(\frac{x}{pq},p\right)$ by $\frac{|C|}{|G|}\int_y^{\frac{x}{q}}\Psi\left(\frac{x}{tq},t\right)\frac{dt}{\log t}$
\end{itemize}
But then, replacing the summations using integrals will lead to some error in equation (2.4). Let us now estimate the errors that can arise due to this change. We denote:
\begin{eqnarray}
    E_1 := \left| \sum_{y < q \leq \sqrt{x}}\left(\sum_{y<p<q;\;\left[\frac{K/\mathbb{Q}}{p}\right]=C} \Psi\left(\frac{x}{pq},p\right) - \frac{|C|}{|G|}\int_y^q\Psi\left(\frac{x}{tq},t\right)\frac{dt}{\log t}\right)\right| \nonumber \\
    \text{and}\;\;\;\;E_2:= \left|\sum_{\sqrt{x}<q\leq \frac{x}{y}}\left(\sum_{y<p<\frac{x}{q};\;\left[\frac{K/\mathbb{Q}}{p}\right]=C} \Psi\left(\frac{x}{pq},p\right) - \frac{|C|}{|G|}\int_y^{\frac{x}{q}}\Psi\left(\frac{x}{tq},t\right)\frac{dt}{\log t}\right)\right| \nonumber
\end{eqnarray}
\underline{\textit{Estimating $E_1$:}} Rewriting $E_1$, we get
\begin{eqnarray}
    E_1 = \left| \sum_{y < q \leq \sqrt{x}}\left(\sum_{y<p<q;\;\left[\frac{K/\mathbb{Q}}{p}\right]=C}\;\; \sum_{n \leq \frac{x}{pq};\;P_1(n)\leq p} 1 - \frac{|C|}{|G|}\int_y^q \left[\sum_{n \leq \frac{x}{pq};\;P_1(n)\leq p} 1\right]\frac{dt}{\log t}\right)\right| \nonumber \\
    \leq \sum_{y < q \leq \sqrt{x}}\;\; \sum_{n \leq \frac{x}{yq}} \left|\sum_{\max(y,P_1(n))\leq p \leq \min\left(q,\frac{x}{nq}\right);\;\left[\frac{K/\mathbb{Q}}{p}\right]=C} 1-\frac{|C|}{|G|}\int_{\max(y,P_1(n))}^{\min\left(q,\frac{x}{nq}\right)}\frac{dt}{\log t}\right| \nonumber
\end{eqnarray}
Now, using Theorem 1.2, we get that 
\begin{eqnarray}
    E_1 << \sum_{y < q \leq \sqrt{x}}\;\; \sum_{n \leq \frac{x}{yq}} \frac{x}{nq\exp \left\{\sqrt{\frac{\log(x/nq)}{n_K}}\right\}} << \frac{x}{\exp\{\sqrt{\log y}\}} \sum_{y<q\leq \sqrt{x}}\frac{1}{q}\sum_{n\leq \frac{x}{yq}}\frac{1}{n} \nonumber \\
    <<\frac{x\log x}{\exp\{\sqrt{\log y}\}} \sum_{y<q\leq \sqrt{x}}\frac{1}{q}<<\frac{x\log x\log\log x}{\exp\{\sqrt{\log y}\}}
\end{eqnarray}
\underline{\textit{Estimating $E_2$:}} Rewriting $E_2$, we get
\begin{eqnarray}
    \left|\sum_{\sqrt{x}<q\leq \frac{x}{y}}\left(\sum_{y<p<\frac{x}{q};\;\left[\frac{K/\mathbb{Q}}{p}\right]=C} \;\;\sum_{n \leq \frac{x}{pq};\;P_1(n)\leq p} 1 - \frac{|C|}{|G|}\int_y^{\frac{x}{q}}\left[\sum_{n \leq \frac{x}{pq};\;P_1(n)\leq p} 1\right]\frac{dt}{\log t}\right)\right| \nonumber \\
    \leq  \sum_{\sqrt{x} < q \leq \frac{x}{y}}\;\; \sum_{n \leq \frac{x}{yq}} \left|\sum_{\max(y,P_1(n))\leq p \leq \frac{x}{nq};\;\left[\frac{K/\mathbb{Q}}{p}\right]=C} 1-\frac{|C|}{|G|}\int_{\max(y,P_1(n))}^{\frac{x}{nq}}\frac{dt}{\log t}\right| \nonumber
\end{eqnarray}
Again using Theorem 1.2, we get that
\begin{eqnarray}
    E_2 <<  \sum_{\sqrt{x} < q \leq \frac{x}{y}}\;\; \sum_{n \leq \frac{x}{yq}} \frac{x}{nq\exp \left\{\sqrt{\frac{\log(x/nq)}{n_K}}\right\}} << \frac{x}{\exp\{\sqrt{\log y}\}} \sum_{\sqrt{x}<q\leq \frac{x}{y}}\frac{1}{q}\sum_{n\leq \frac{x}{yq}}\frac{1}{n}  \nonumber \\
    << \frac{x\log x}{\exp\{\sqrt{\log y}\}} \sum_{\sqrt{x}<q\leq \frac{x}{y}}\frac{1}{q}<< \frac{x\log x \log \log x}{\exp\{\sqrt{\log y}\}} 
\end{eqnarray}
Therefore, from (2.4),(2.5) and (2.6), we get that
\begin{eqnarray}
     \sum_{y<p \leq \sqrt{x};\;\left[\frac{K/\mathbb{Q}}{p}\right]=C} |S_2(x,p)| = \frac{|C|}{|G|}\left[ \sum_{y < q \leq \sqrt{x}}\int_y^q\Psi\left(\frac{x}{tq},t\right)\frac{dt}{\log t} + \sum_{\sqrt{x}<q\leq \frac{x}{y}}\;\;\int_y^{\frac{x}{q}}\Psi\left(\frac{x}{tq},t\right)\frac{dt}{\log t} \right] \nonumber \\
    + O\left(\frac{x\log x\log \log x}{\exp\{\sqrt{\log y}\}}\right) \hspace{1.5cm }
\end{eqnarray}
It is now important to evaluate these integrals to get an estimate for the sum of $|S_2(x,p)|$, and therefore, prove the theorem. In pursuit of that, we observe, again by Theorem 1.4, that
\begin{eqnarray}
    \sum_{p \leq \sqrt{x}} |S_2(x,p)| = [x]-\pi(x)+O\left(\frac{x}{\exp\{\frac{1}{2}\sqrt{\log x \log \log x}\}}\right)  = x+O\left(\frac{x}{\log x}\right)
\end{eqnarray}
Using Theorem 1.5 and (2.8), we then get that
\begin{eqnarray}
    x+O\left(\frac{x\log y}{\log x}\right) = \sum_{y<p\leq \sqrt{x}} |S_2(p,x)| = \sum_{y < q \leq \sqrt{x}}\;\;\sum_{y<p<q} \Psi\left(\frac{x}{pq},p\right) + \sum_{\sqrt{x}<q\leq \frac{x}{y}}\;\;\sum_{y<p\leq \frac{x}{q}} \Psi\left(\frac{x}{pq},p\right)   
\end{eqnarray}
Note that the summation split on the right of (2.9) is exactly the same as that of (2.4) without the conjugacy class condition. Using the strong form of the Prime Number Theorem to estimate the error, a similar replacement of the inner sums in (2.9) with the integrals as done above would then yield the following from (2.9):
\begin{eqnarray}
     \sum_{y < q \leq \sqrt{x}}\int_y^q\Psi\left(\frac{x}{tq},t\right)\frac{dt}{\log t} + \sum_{\sqrt{x}<q\leq \frac{x}{y}}\;\;\int_y^{\frac{x}{q}}\Psi\left(\frac{x}{tq},t\right)\frac{dt}{\log t} = x+O\left(\frac{x\log y}{\log x}\right) + O\left(\frac{x\log x\log \log x}{\exp\{\sqrt{\log y}\}}\right) \hspace{0.5cm}
\end{eqnarray}
Therefore, from (2.3), (2.7) and (2.9), and taking $y = \exp{(2\log\log x)^2}$ we get that
\begin{eqnarray}
    |S_2(x,C)| = \frac{|C|}{|G|}.x + O\left(\frac{x(\log\log x)^2}{\log x}\right) 
\end{eqnarray}
But then,
\begin{eqnarray}
    N_2(x;K/\mathbb{Q},C) = |S_2(x,C)| 
\end{eqnarray}
Hence, finally adjoining (2.11) and (2.12), we get our desired result:
\begin{eqnarray}
      N_2(x;K/\mathbb{Q},C)=\sum_{n \leq x;\;\left[\frac{K/\mathbb{Q}}{P_2(n)}\right]=C } 1 = \frac{|C|}{|G|}.x + O\left(\frac{x(\log \log x)^2}{\log x}\right) 
\end{eqnarray}
We will now prove the following lemma that we will use in the theorems that will follow towards our main result. \\ \\
\textbf{Lemma 2.2:} \textit{For a finite Galois extension $K$ of $\mathbb{Q}$ and a conjugacy class $C$ of the corresponding Galois group, we have}
\begin{eqnarray}
     \sum_{n \leq x;\;\left[\frac{K/\mathbb{Q}}{p_1(n)}\right]=C} \mu(n) << \frac{x\log x}{\exp\{(\log x)^{\frac{1}{3}}\}} \nonumber
\end{eqnarray}
\textbf{Proof:} Let $f$ be an arithmetic function defined as 
\begin{eqnarray}
    f(m)=\begin{cases}
    1,\;\text{if $m$ is a prime and $\left[\frac{K/\mathbb{Q}}{m}\right]=C$ } \\
    0, \;\text{otherwise}
    \end{cases} \nonumber
\end{eqnarray}
Therefore, we have 
\begin{eqnarray}
   \sum_{n \leq x;\;\left[\frac{K/\mathbb{Q}}{p_1(n)}\right]=C} \mu(n)=  \sum_{n \leq x} \mu(n)f(p_1(n)) 
\end{eqnarray}
Using Möbius inversion on Alladi Duality Formula for $k=1$(here we use the second identity in Theorem 1.3) and then the hyperbola method for splitting the summation, we get from above
\begin{eqnarray}
    \sum_{n \leq x} \mu(n)f(p_1(n)) = -\sum_{n \leq x}\sum_{d|n}\mu\left(\frac{n}{d}\right)f(P_1(d)) \hspace{5.6cm}\nonumber \\
    =  -\sum_{ m \leq \sqrt{x}} \mu(m)\sum_{d \leq \frac{x}{m}}f(P_1(d))-\sum_{d \leq \sqrt{x}}f(P_1(d))\sum_{\sqrt{x}<m<\frac{x}{d}}\mu(m) 
\end{eqnarray}
We now recall Theorem 2 of [4]. With the same hypothesis as of the lemma, we have that
\begin{eqnarray}
    \sum_{2 \leq n \leq x;\;\left[\frac{K/\mathbb{Q}}{P_1(n)}\right]=C} 1 = \frac{|C|}{|G|}.x+O\left(\frac{x}{\exp\{(\log x)^{\frac{1}{3}}\}}\right) 
\end{eqnarray}
Further, we recall the following two consequences of the strong form of Prime Number Theorem:
\begin{eqnarray}
     \sum_{1 \leq n \leq x} \mu(n) << \frac{x}{\exp\{\sqrt{\log x}\} } \;\;\;\;\text{and}\;\;\;\;
 \sum_{1 \leq n \leq x} \frac{\mu(n)}{n} << \frac{1}{\exp\{\sqrt{\log x}\}}
\end{eqnarray}
Now, continuing with the first term in the RHS of (2.15) and using (2.16) and (2.17), we have
\begin{eqnarray}
    -\sum_{ m \leq \sqrt{x}} \mu(m)\sum_{d \leq \frac{x}{m}}f(P_1(d)) = -\sum_{m \leq \sqrt{x}}\mu(m)\left[\frac{|C|}{|G|}.\frac{x}{m}+O\left(\frac{x}{m\exp\{(\log\frac{x}{m})^{\frac{1}{3}}\}}\right) \right] \hspace{0.1cm}\nonumber \\
    = -\frac{|C|}{|G|}.O\left(\frac{x}{\exp\{(\log x)^{\frac{1}{2}}\}}\right) + O\left(\frac{x\log x}{\exp\{(\log x)^{\frac{1}{3}}\}}\right) \nonumber \\
    << \frac{x\log x}{\exp\{(\log x)^{\frac{1}{3}}\}} \hspace{6cm}
\end{eqnarray}
Now, for the second term in RHS of (2.15), we have using (2.17) that
\begin{eqnarray}
    -\sum_{d \leq \sqrt{x}}f(P_1(d))\sum_{\sqrt{x}<m<\frac{x}{d}}\mu(m) << \sum_{d \leq \sqrt{x}} \frac{x}{d\exp\{\sqrt{\log \frac{x}{d}}\}}  << \frac{x\log x}{\exp\{\sqrt{\frac{\log x}{2}}\}}
\end{eqnarray}
Therefore, combining equations (2.14), (2.15), (2.18) and (2.19), we get our desired result. \qed \\ \\
\textbf{Remark:} Dawsey's derivation of her main result does not use the estiamte in \textit{Lemma 2.2}. She utilizes an alternate derivation due to [Alladi 1] using sums involving $\mu(n)\log n$. But, that approach applied to $P_2(n)$ presents difficulties. Hence, we have established \textit{Lemma 2.2} following Alladi's original derivation.\\ \\ 
As a major step towards our final result, and also a consequence of both Theorem 2.1 and Lemma 2.2, we have our next theorem. \\ \\
\textbf{Theorem 2.3}  \textit{For a finite Galois extension $K$ of $\mathbb{Q}$ and a conjugacy class $C$ of the corresponding Galois group, we have}
\begin{eqnarray}
   \mathcal{S}_{\omega}(x;K/\mathbb{Q},C):= \sum_{n \leq x;\;\left[\frac{K/\mathbb{Q}}{p_1(n)}\right]=C} \mu(n)\omega(n)<< \frac{x(\log\log x)^4}{\log x} \nonumber
\end{eqnarray}
\textbf{Proof:} 
Firstly, we have 
\begin{eqnarray}
    \mathcal{S}_{\omega}(x;K/\mathbb{Q},C) = \sum_{n \leq x}\mu(n)\omega(n)f(p_1(n))
\end{eqnarray}
Note that using Möbius inversion in (1.3), we get that for positive integers $n$,
\begin{eqnarray}
    \mu(n)(\omega(n)-1)f(p_1(n)) = \sum_{d|n}\mu\left(\frac{n}{d}\right)f(P_2(d)) 
\end{eqnarray}
Hence, from (2.20) and (2.21), we get that
\begin{eqnarray}
      \mathcal{S}_{\omega}(x;K/\mathbb{Q},C) = \sum_{n \leq x}\sum_{d|n}\mu\left(\frac{n}{d}\right)f(P_2(d)) + \sum_{n \leq x}\mu(n)f(p_1(n)) \nonumber
\end{eqnarray}
 Therefore, using Lemma 2.2 in the above equation yields
\begin{eqnarray}
          \mathcal{S}_{\omega}(x;K/\mathbb{Q},C) = \sum_{n \leq x}\sum_{d|n}\mu\left(\frac{n}{d}\right)f(P_2(d)) + O\left(\frac{x\log x}{\exp\{(\log x)^{\frac{1}{3}}\}}\right)
\end{eqnarray}
Applying the hyperbola method in the first term of RHS of (2.22), we get
\begin{eqnarray}
    \sum_{n\leq x}\sum_{d|n}\mu\left(\frac{n}{d}\right)f(P_2(d)) = \sum_{m \leq T}\mu(m)\sum_{d \leq \frac{x}{m}}f(P_2(d))+\sum_{d \leq \frac{x}{T}}f(P_2(d))\sum_{T<m\leq \frac{x}{d}} \mu(m) 
\end{eqnarray}
Using the first bound of (2.17) in the second term on the RHS of (2.23), we get that
\begin{eqnarray}
    \sum_{d \leq \frac{x}{T}}f(P_2(d))\sum_{T<m\leq \frac{x}{d}} \mu(m) << \sum_{d \leq \frac{x}{T}}f(P_2(d))\frac{x}{d\exp\left\{\sqrt{\log \left(\frac{x}{d}\right)}\right\}}\hspace{1.2cm} \nonumber\\ << \frac{x}{\exp\{\sqrt{\log T}\}}\sum_{d \leq \frac{x}{T}}\frac{1}{d} <<\frac{x\log x}{\exp\{\sqrt{\log T}\}}
\end{eqnarray}
Now, for the other summation in RHS of (2.23), we get using Theorem 2.1 that
\begin{eqnarray}
    \sum_{m \leq T}\mu(m)\sum_{d \leq \frac{x}{m}}f(P_2(d))= \frac{|C|}{|G|}x\sum_{m \leq T}\frac{\mu(m)}{m} + O\left(x(\log\log x)^2\sum_{m\leq T}\frac{\mu(m)}{m\log\left(\frac{x}{m}\right)}\right)
\end{eqnarray}
Using the second bound of (2.17) in (2.25), we then get that
\begin{eqnarray}
    \sum_{d \leq \frac{x}{T}}f(P_2(d))\sum_{T<m\leq \frac{x}{d}} \mu(m) << \frac{|C|}{|G|}\frac{x}{\exp\{\sqrt{\log T}\}}+\frac{x\log T(\log\log x)^2}{\log \left(\frac{x}{T}\right)}
\end{eqnarray}
Therefore, (2.22),(2.23),(2.24) and (2.26) together yield
\begin{eqnarray}
    \mathcal{S}_{\omega}(x;K/\mathbb{Q},C) << \frac{x\log x}{\exp\{\sqrt{\log T}\}} + \frac{|C|}{|G|}\frac{x}{\exp\{\sqrt{\log T}\}}+\frac{x\log T(\log\log x)^2}{\log \left(\frac{x}{T}\right)} + \frac{x\log x}{\exp\{(\log x)^{\frac{1}{3}}\}}
\end{eqnarray}
Choosing $T=\exp\{(2\log\log x)^2\}$, we finally get our desired reuslt:
\begin{eqnarray}
    \mathcal{S}_{\omega}(x;K/\mathbb{Q},C)=\sum_{n \leq x;\;\left[\frac{K/\mathbb{Q}}{p_1(n)}\right]=C} \mu(n)\omega(n) << \frac{x(\log\log x)^4}{\log x} \nonumber
\end{eqnarray} \qed \\ \\ 
\textbf{Theorem 2.4} \textit{For a finite Galois extension $K$ of $\mathbb{Q}$ and a conjugacy class $C$ of the corresponding Galois group, we have}
\begin{eqnarray}
    \sum_{n \leq x;\;\left[\frac{K/\mathbb{Q}}{p_1(n)}\right]=C} \mu(n)\omega(n) \left\{\frac{x}{n}\right\} << \frac{x(\log\log x)^{\frac{5}{2}}}{\sqrt{\log x}} \nonumber
\end{eqnarray} 
where $\{.\}$ denotes the fractional part. \\ \\
\textbf{Proof:} First we decompose the summation in our desired result as follows:
\begin{eqnarray}
    \sum_{n \leq x;\;\left[\frac{K/\mathbb{Q}}{p_1(n)}\right]=C} \mu(n)\omega(n) \left\{\frac{x}{n}\right\} = \sum_{n \leq T;\;\left[\frac{K/\mathbb{Q}}{p_1(n)}\right]=C} \mu(n)\omega(n) \left\{\frac{x}{n}\right\}+\sum_{T<n \leq x;\;\left[\frac{K/\mathbb{Q}}{p_1(n)}\right]=C} \mu(n)\omega(n) \left\{\frac{x}{n}\right\}
\end{eqnarray}
where $T$ is an appropriate function of $x$ to be chosen at the end of the proof.
Looking at the first term of the RHS of (2.28), we get that
\begin{eqnarray}
    \sum_{n \leq T;\;\left[\frac{K/\mathbb{Q}}{p_1(n)}\right]=C} \mu(n)\omega(n) \left\{\frac{x}{n}\right\}<<\sum_{n \leq T}\omega(n) << T\log\log T 
\end{eqnarray}
Now, using Theorem 2.3, we have from above that
\begin{eqnarray}
    \sum_{T<n \leq x;\;\left[\frac{K/\mathbb{Q}}{p_1(n)}\right]=C} \mu(n)\omega(n) \left\{\frac{x}{n}\right\} = \sum_{T < n \leq x} [\mathcal{S}_{\omega}(n,K/\mathbb{Q},C)-\mathcal{S}_{\omega}(n-1,K/\mathbb{Q},C)]\left\{\frac{x}{n}\right\}  \nonumber\\
    = \sum_{T<n\leq x} \mathcal{S}_{\omega}(n,K/\mathbb{Q},C)\left(\left\{\frac{x}{n}\right\}-\left\{\frac{x}{n+1}\right\}\right) \hspace{1.2cm} \nonumber \\
    << \frac{x(\log \log x)^4}{\log x}V_{\{\}}\left[1,\frac{x}{T}\right]
    << \frac{x(\log \log x)^4}{\log x}.\frac{x}{T} \hspace{0.8cm}
\end{eqnarray}
where $V_{\{\}}[a,b]$ is the bounded variation of $\{x\}$ on the interval $[a,b]$. 
Combining (2.28)-(2.30), we then get
\begin{eqnarray}
    \sum_{n \leq x;\;\left[\frac{K/\mathbb{Q}}{p_1(n)}\right]=C} \mu(n)\omega(n) \left\{\frac{x}{n}\right\} << T\log\log T+\frac{x(\log \log x)^4}{\log x}.\frac{x}{T} \nonumber
\end{eqnarray}
Choosing $T = \frac{x(\log\log x)^{\frac{3}{2}}}{\sqrt{\log x}}$ then gives us our desired bound:
\begin{eqnarray}
       \sum_{n \leq x;\;\left[\frac{K/\mathbb{Q}}{p_1(n)}\right]=C} \mu(n)\omega(n) \left\{\frac{x}{n}\right\} << \frac{x(\log\log x)^{\frac{5}{2}}}{\sqrt{\log x}} \nonumber
\end{eqnarray} \qed \\ \\
We now prove our final result of this paper. \\ \\
\textbf{Theorem 2.5:} \textit{For a finite Galois extension $K$ of $\mathbb{Q}$ and a conjugacy class $C$ of the corresponding Galois group, we have}
\begin{eqnarray}
    \sum_{n \leq x;\;\left[\frac{K/\mathbb{Q}}{p_1(n)}\right]=C} \frac{\mu(n)\omega(n)}{n}<<\frac{(\log\log x)^{\frac{5}{2}}}{\sqrt{\log x}} \nonumber
\end{eqnarray}
\textit{More precisely, taking $x \rightarrow \infty$, we get}
\begin{eqnarray}
    \sum_{1 \leq n;\;\left[\frac{K/\mathbb{Q}}{p_1(n)}\right]=C} \frac{\mu(n)\omega(n)}{n}=\sum_{2 \leq n;\;\left[\frac{K/\mathbb{Q}}{p_1(n)}\right]=C} \frac{\mu(n)\omega(n)}{n}=0 \nonumber
\end{eqnarray}
To prove Theorem 2.5, we first prove the following lemma: \\ \\
\textbf{Lemma 2.6:} \textit{With the hypothesis as in Theorem 2.4, we have}
\begin{eqnarray}
    \sum_{n\leq x\;\left[\frac{K/\mathbb{Q}}{p_1(n)}\right]=C}\mu(n)\omega(n)\left[\frac{x}{n}\right] << \frac{x(\log\log x)^2}{\log x} \nonumber
\end{eqnarray}
\textbf{Proof:} Using the function $f$ already defined in the proof of Lemma 2.2, we get that
\begin{eqnarray}
 \sum_{n\leq x\;\left[\frac{K/\mathbb{Q}}{p_1(n)}\right]=C}\mu(n)\omega(n)\left[\frac{x}{n}\right] =  \sum_{n\leq x}\mu(n)\omega(n)f(p_1(n))\left[\frac{x}{n}\right]
\end{eqnarray}
From Theorem 2.1 and (1.3), we get
\begin{eqnarray}
    \sum_{n\leq x}\mu(n)(\omega(n)-1)f(p_1(n))\left[\frac{x}{n}\right] = \sum_{n \leq x}\sum_{d|n}\mu(d)(\omega(d)-1)f(p_1(d))\hspace{2.1cm}\nonumber \\
    =\sum_{n\leq x}f(P_2(n)) =   \frac{|C|}{|G|}.x + O\left(\frac{x(\log \log x)^2}{\log x}\right)
\end{eqnarray}
Further, using (2.16) and Alladi's Duality formula for $k=1$ case, we get that
\begin{eqnarray}
    \sum_{n \leq x}\mu(n)f(p_1(n))\left[\frac{x}{n}\right] = \sum_{n \leq x}\sum_{d|n} \mu(d)f(p_1(d))=-\sum_{n \leq x}f(P_1(n)) = -\frac{|C|}{|G|}.x + O\left(x\exp\{-(\log x)^{\frac{1}{3}}\}\right) \hspace{0.5cm}
\end{eqnarray}
Combining (2.31)-(2.33), we get
\begin{eqnarray}
    \sum_{n\leq x\;\left[\frac{K/\mathbb{Q}}{p_1(n)}\right]=C}\mu(n)\omega(n)\left[\frac{x}{n}\right] =O\left(\frac{x(\log \log x)^2}{\log x}\right)+ O\left(x\exp\{-k(\log x)^{\frac{1}{3}}\}\right)<< \frac{x(\log \log x)^2}{\log x} \nonumber 
\end{eqnarray}
which proves the lemma. \qed \\ \\
\textbf{Proof of Theorem 2.5:} 
Combining Theorem 2.4 and Lemma 2.6, we get
\begin{eqnarray}
     \sum_{n \leq x;\;\left[\frac{K/\mathbb{Q}}{p_1(n)}\right]=C} \mu(n)\omega(n) \left\{\frac{x}{n}\right\}+\sum_{n\leq x\;\left[\frac{K/\mathbb{Q}}{p_1(n)}\right]=C}\mu(n)\omega(n)\left[\frac{x}{n}\right] << \frac{x(\log\log x)^{\frac{5}{2}}}{\sqrt{\log x}} +\frac{x(\log \log x)^2}{\log x} \nonumber \\
     \implies x \sum_{n \leq x;\;\left[\frac{K/\mathbb{Q}}{p_1(n)}\right]=C} \frac{\mu(n)\omega(n)}{n} << \frac{x(\log\log x)^{\frac{5}{2}}}{\sqrt{\log x}} \implies \sum_{n \leq x;\;\left[\frac{K/\mathbb{Q}}{p_1(n)}\right]=C} \frac{\mu(n)\omega(n)}{n} << \frac{(\log\log x)^{\frac{5}{2}}}{\sqrt{\log x}} 
\end{eqnarray}
This proves the first part of the theorem. Hence, taking $x \rightarrow \infty$, we finally get that
\begin{eqnarray}
    \sum_{1 \leq n;\;\left[\frac{K/\mathbb{Q}}{p_1(n)}\right]=C} \frac{\mu(n)\omega(n)}{n} = 0 \nonumber
\end{eqnarray}
But then, we note that $\omega(1)=0$. Therefore,
\begin{eqnarray}
    \sum_{2 \leq n;\;\left[\frac{K/\mathbb{Q}}{p_1(n)}\right]=C} \frac{\mu(n)\omega(n)}{n} = 0 \nonumber
\end{eqnarray}\qed \\ \\
Thus, we have established our main result that for any arbitrary finite Galois extension $K$ of the field of rationals $\mathbb{Q}$ and a conjugacy class $C$ of the Galois group $G=Gal(K/\mathbb{Q})$, 
\begin{eqnarray}
    \sum_{2 \leq n;\;\left[\frac{K/\mathbb{Q}}{p_1(n)}\right]=C} \frac{\mu(n)\omega(n)}{n} = 0 
\end{eqnarray}

\section{Examples and Illustrations}
To understand the actual implication of results, there is always a need of looking at proper examples relating to them. So in this section, we focus on looking at a few examples that will illustrate the results proven in the previous section.
\begin{itemize}
    \item[1.] Let us consider $K=\mathbb{Q}(\zeta_k)$, i.e. the $k^{th}$ cyclotomic field extension of $\mathbb{Q}$. The Galois group of the corresponding extension is $\mathbb{Z}_k^*$. Since the Galois group is abelian, every conjugacy class in the group is a single element. So, let $C$ denote a conjugacy class in the Galois group. More precisely, let $C$ denote the map $\sigma_l$ in the Galois group, such that $\sigma_l(\zeta_k)=\zeta_k^l$. Now, in this setting, we look at the condition which we use to slice the summation in our main result. That is, we look at the numbers $n$, such that their smallest prime factor $p_1(n)$ is unramified in $K$ and the corresponding Artin Symbol $\left[\frac{K/\mathbb{Q}}{p_1(n)}\right]=\sigma_l$. Firstly, for $p_1(n)$ to be unramified in $K$, $p_1(n)$ does not divide $ k$. Now by definition, we know that the Artin Symbol is the Frobenius element, and in this example, it is the map $\sigma_{p_1(n)}$ given by $\sigma_{p_1(n)}(\zeta_k)=\zeta_k^{p_1(n)}$. To satisfy the condition, we must have that $\sigma_{p_1(n)}=\sigma_l$, i.e. $p_1(n)\equiv l\;(mod \;k)$. This is exactly the condition on integers for the main result of [2]. Therefore, this example is the case for Alladi-Johnson's result and therefore, we know that Theorem 2.5 is true in this setting. Note that since $p_1(n)$ is unramified and $p_1(n) \equiv l\;(mod\;k)$, we automatically have that $(l,k)=1$.
    \item[2.] We have also worked out the case when $K$ is the splitting field of the polynomial $x^3+x+1$. Of course, the extension is finite and the Galois group of the extension is $S_3$. There are 3 different conjugacy classes of $S_3$. Precisely, they are comprised of $1-$cycles, $2-$cycles and $3-$cycles respectively. According to Theorem 2.5, we have to fix a conjugacy class, and hence we have three different cases. 
    \begin{itemize}
        \item[i)] Let $C$ be the conjugacy class containing the $3-$cycles. So for the sum, we only need to consider all those integers for which, the Artin symbol corresponding to their smallest prime factors is a $3-$cycle. Now, speaking algebraically, if the Artin Symbol of $p_1(n)$ is a $3-$cycle, then the polynomial $x^3+x+1$ is irreducible in the residue field $\mathbb{Z}/p_1(n)\mathbb{Z}$. 
        \item[ii)] Let $C$ now be the conjugacy class containing $2-$cycles. In this case, we look at all those integers for which the corresponding Artin Symbol of their smallest prime factors is a $2-$cycle. Those are all those integers $n$, for which the polynomial $x^3+x+1$ has exactly one root in the residue field  $\mathbb{Z}/p_1(n)\mathbb{Z}$. 
        \item[iii)] The final case is that of the chosen conjugacy class is the singleton set, consisting of the $1-$cycle, i.e. the identity. In this case, we take the sum over all those integers $n$, for which the polynomial completely splits in the  residue field  $\mathbb{Z}/p_1(n)\mathbb{Z}$. 
\end{itemize}
We calculated the sum on the LHS of (2.35) for each of these case, fixing the conjugacy classes separately, for numbers up to 80,000. The table below show the results that we have received using SageMath:
\begin{center}
    \begin{tabular}{|c|c|c|c|}
    \hline
        Table for Sum estimates & $n \leq 20000$ & $n \leq 40000$ & $n \leq 80000$ \\
        \hline 
        $3-$cycles & 0.25 &  0.254 & 0.265 \\
        \hline 
        $2-$cycles & 0.188 & 0.237 & 0.279 \\
        \hline 
        $1-$cycle or identity & -0.026& -0.008& 0.009 \\
        \hline
    \end{tabular}
\end{center}
\end{itemize}
\textbf{Note:} We note here that the convergence of the sum is really slow since the estimated error is given by $\frac{(\log \log x)^{5/2}}{\sqrt{\log x}}$. Given that, the results noted in the above table are quite interesting and show that indeed the sum is converging to $0$. \\ \\
\textbf{Remark:} If we look at the absolute sum, then the Moebius function does not play any role, and hence, the sum $\sum_{n \leq x;\;\left[\frac{K/\mathbb{Q}}{p_1(n)}\right]=C} \frac{\omega(n)}{n}$ is quite large. So, we owe it to $\mu(n)$ for the cancellations that occur, and hence the convergence of the sum.
\section{The Integers that have not been Considered}
We look at the conclusion of \textit{Theorem 2.5} and see that the infinite sum on the left hand side, takes into consideration only those integers that satisfy the condition of the Artin Symbol corresponding to the fixed conjugacy class. Of course, that set of integers is only but a proper subset of the set of all the positive integers, i.e. there exists integers that are not considered in the sum. This is automatic, since once fixing the conjugacy class would force us to consider only those integers that would yield the corresponding Artin Symbol of their smallest prime factors to be the fixed conjugacy class, and the rest are left out. Precisely, if we fix a conjugacy class $C$ of the said Galois group, then for every integer $n$ such that $\left[\frac{K/\mathbb{Q}}{p_1(n)}\right] \ne C$, $n$  is not involved in the sum. Further, there are integers that are not considered in any of the sums, for different conjugacy classes. Those integers are precisely the ones, whose smallest prime factors are \textit{ramified in $K$}.  \\ \\
Now, let us look at the following \textit{Theorem 4} in \S 2 in [2]. \\ \\
\textbf{Theorem:}(Alladi-Johnson) \textit{We have that}
\begin{eqnarray}
    \sum_{n \geq 2} \frac{\mu(n)\omega(n)}{n}=0 
\end{eqnarray}
Note that, the summand in (4.1) is the exactly the same as our Theorem 2.5. The only difference that lies is that we have summed under the condition of the Artin Symbol, which has led us to take partial infinite sums of (4.1). The most surprising and interesting fact is that even though we are considering only partial sums under a particular condition, we are still getting the same value, i.e. 0.\\ \\
The contribution in the sum in (4.1) from all the integers involved in either of the sums, for different conjugacy classes, are all 0, by Theorem 2.5. We then try to understand the contribution of all the integers whose smallest prime factors are ramified in $K$. The follow \textit{Theorem 11} in \S7 in [2] will give us what we need:\\ \\
\textbf{Theorem:} \textit{Let $p$ be an arbitrary but fixed prime. Then}
\begin{eqnarray}
    \sum_{n\geq 2;\;p_1(n)=p} \frac{\mu(n)\omega(n)}{n}=0
\end{eqnarray}
The above theorem holds for any arbitrary prime $p$. Since, there are only finitely primes that are ramified in $K$, we therefore conclude that the contribution from all those finitely many primes is also 0, thus, making sense of the complete sum in (4.1). \\ \\
\textbf{Remark:} The most surprising of the facts is that each of the partial sums that add up to the total sum in (4.1), be it the integers corresponding  to the fixed conjugacy classes or the integers with ramified smallest prime factor, all contribute nothing but 0, thus making the total sum 0. \\ Further, in Theorem 2.5, we are really \textit{adding infinitely many zeros, to get zero}!
\section{A Few Remarks Corresponding to the General $f$}
In [2], the authors state the following Theorem 13 in \S8: \\ \\
\textbf{Theorem:} \textit{If $f$ is a bounded function on the prime such that} 
\begin{eqnarray}
    \sum_{2 \leq n \leq x}f(P(n)) \sim \kappa x \nonumber
\end{eqnarray}
\textit{and}
\begin{eqnarray}
     \sum_{2 \leq n \leq x}f(P_2(n)) \sim \kappa x \nonumber
\end{eqnarray}
\textit{for some constant $\kappa$, then}
\begin{eqnarray}
    \sum_{n=2}^{\infty} \frac{\mu(n)\omega(n)f(n)}{n}=0 \nonumber
\end{eqnarray}
Now, looking at our case, we defined the function $f$ in Lemma 2.2, which is bounded and non-zero on primes. We observe that \textit{Theorem 2} of [4] gives us an average of $f(P(n))$ and our \textit{Theorem 2.1} gives us an average of $f(P_2(n))$. Coincidentally, or rather not, we find that the averages in both the cases are same, that is $|C|/|G|$. Therefore, our case also illustrates the above Theorem proved by Alladi and Johnson. Precisely, our Theorem 2.5 is exactly the conclusion of the above Theorem by A-J. \\ \\
\textbf{Note:} Our result does not stand insignificant in the face of the above theorem, since we had to prove the average of $f(P_2(n))$ is $|C|/|G|$, which is Theorem 2.1 in \S2. The rest of the proofs to our main goal follows quite a similar approach as of the above Theorem by A-J. Also, our choice of $f$ correctly fits in the context.
\newpage
\section{Arithmetic Density Versions}
The generalizations of (1.4) to algebraic number fields by various authors starting with Dawsey [4] was motivated by rewriting (1.4) as
\begin{eqnarray}
-\sum_{n\ge 2, \, p(n)\equiv\ell(mod\,k)}\frac{\mu(n)}{n}=\frac{1}{\phi(k)},    
\end{eqnarray}
and interpreting this as an arithmetic density result. As we have already seen in (1.5), Dawsey's result actually gives a new formula for the Chebotarev Density, when rewritten as follows:
\begin{eqnarray}
      -\lim_{x \to \infty}\sum_{\substack{2 \leq n \leq x \\ \left[\frac{K/\mathbb{Q}}{p_1(n)}\right] = C}} \frac{\mu(n)}{n} = \frac{|C|}{|G|}.
\end{eqnarray}
In the light of (6.2), we can also rewrite our main result (2.35) as
\begin{eqnarray}
      \sum_{2 \leq n;\;\left[\frac{K/\mathbb{Q}}{p_1(n)}\right]=C} \frac{\mu(n)(\omega(n)-1)}{n} =   \sum_{2 \leq n;\;\left[\frac{K/\mathbb{Q}}{p_1(n)}\right]=C} \frac{\mu(n)\omega(n)}{n} -   \sum_{2 \leq n;\;\left[\frac{K/\mathbb{Q}}{p_1(n)}\right]=C} \frac{\mu(n)}{n} = \frac{|C|}{|G|}.
\end{eqnarray}
So, (6.3) is the density version of (2.35) under the specified conditions in \S 2. We also note that, the LHS of (6.3) gives us a new formula for the Chebotarev Density.

\section{Further Results Obtained}
In 2019, Sweeting and Woo [11] proved algebraic analogues of Alladi's 1977 result for general number field extensions. A major result in that paper is the duality lemma for general number fields. Just like in Alladi's 1977 paper, the duality lemma for general number fields has paved a way for the general extension of Dawsey's result. \\ \\
Encouraged by these results, the author of this paper went on to establish a general duality lemma [10] (like that has been proved by Alladi in case of integers in 1977) involving prime ideals of smallest and $k^{th}$ largest norms. This general duality for arbitrary Galois extension of number fields will be stepping stone in proving analogues of results in [2] in arbitary extensions, and not only in rational extensions, which is already proved above. \\ \\
Ideas for future does not end there as the author and his Doctoral Advisor, Dr. K. Alladi, have recently obtained results on higher order dualities [3] which use the general duality lemma proved in [1]. In light of the general duality in number fields, we aim to prove similar higher order duality analogues for arbitrary number field extensions, and we have made significant advances in that direction. This will lead us to have a beautiful and bigger picture of the importance of duality, not only in just the set of integers, but also in general algebraic structures. \\ \\
\textbf{\textit{Acknowledgement:}} I am indebted to my Doctoral Advisor, Dr. K. Alladi for his support and guidance throughout this work, and to one of my professors, Dr. J. Booher for constantly helping with several algebraic explanations and examples required for the work.
\newpage
\section{References}
\begin{itemize}
    \item[1)] K. ALLADI, ``Duality between prime factors and an application to the prime number theorem for arithmetic progressions", \textit{J. Num. Th.}, \textbf{9} (1977), 436-451.
    \item[2)] K. ALLADI and J. JOHNSON, ``Duality between prime factors and the prime number theorem for arithmetic progressions - II", (\textit{submitted pre-print}), 2024. \href{https://nam10.safelinks.protection.outlook.com/?url=http%3A%2F%2Farxiv.org%2Fabs%2F2410.18259&data=05%7C02%7Csengupta.s%40ufl.edu%7Cbb6e0dfda02940b30dbb08dcf7758568%7C0d4da0f84a314d76ace60a62331e1b84%7C0%7C0%7C638657331341043920%7CUnknown%7CTWFpbGZsb3d8eyJWIjoiMC4wLjAwMDAiLCJQIjoiV2luMzIiLCJBTiI6Ik1haWwiLCJXVCI6Mn0%3D%7C0%7C%7C%7C&sdata=YQy1HVSyFY7D%2F13FeVYxRe5BIEnZ%2BmAsVQuJ4xLqiOY%3D&reserved=0}{arXiv:2410.18259}
    \item[3)] K. ALLADI and S. SENGUPTA, ``Higher order duality between prime factors and primes in arithmetic progressions", (\textit{in preparation}).
    \item[4)] M. L. DAWSEY, ``A new formula for Chebotarev densities", \textit{Res. in Num. Th.}, \textbf{3} (2017), 1-13. 
    \item[5)] N. G. DE BRUIJN, ``The asymptotic behaviour of a function occurring in the theory of primes", \textit{J. Indian Math. Soc. (N.S.)}, \textbf{15} (1951), 25-32.
    \item[6)] N. G. DE BRUIJN, ``On the number of positive integers $\leq x$ and free of prime factors $>y$", \textit{Indag. Math.}, \textbf{13} (1951), 50-60.
    \item[7)] A. HILDEBRAND, ``On the number of positive integers $\leq x$ and free of prime factors $>y$", \textit{J. Num. Th.} \textbf{22} (1986), 289-307.
    \item[8)] J. C. LAGARIAS and A. M. ODLYZKO, ``Effective versions of the Chebotarev density theorem", \textit{Algebraic Number Fields: L-functions and Galois properties (A. Fr\"ohlich, Acad. Press. London)}, \textbf{7} (1977), 409-464.
    \item[9)] P. MOREE, ``An interval result for the number field $\Psi(x,y)$ function", \textit{Manuscripta Mathematica}, \textbf{76} (1992), 437-450.
    \item[10)] S. SENGUPTA, ``Higher order dualities in general number fields", (\textit{in preparation}).
    \item[11)] N. SWEETING and K. WOO, ``Formulas for Chebotarev densities of Galois extensions of number fields", \textit{Res. in Num. Th.}, \textbf{5} (2019), 1-13.
    \item[12)] G. TENENBAUM, ``A rate estimate in Billingsley's theorem for the size distribution of large prime factors", \textit{Quart. J. Math.} \textbf{51} (2000), 385-403.
    \item[13)] G. TENENBAUM, ``Introduction to analytic and probabilistic number theory", \textit{Amer. Math. Soc., Providence RI}, \textbf{163} (2015) 
    \item[14)] N. TSCHEBOTAREFF, ``Die Bestimmung der Dichtigkeit einer Menge von Primzahlen, welche zu einer gegebenen Substitutionsklasse geh\"oren", \textit{Math. Ann.} \textbf{95} (1926), 191-228 
\end{itemize}

\end{document}